\documentclass{amsart}
\usepackage{graphicx,amssymb,epsfig}

\usepackage{amssymb,amsmath}
\usepackage[all]{xy}

\newtheorem{theorem}{Theorem}[section]
\newtheorem{proposition}[theorem]{Proposition}
\newtheorem{lemma}[theorem]{Lemma}
\newtheorem{corollary}[theorem]{Corollary}
\newtheorem{remark}[theorem]{Remark}
\newtheorem{example}[theorem]{Example}

\newtheorem{definition}[theorem]{Definition}
\newtheorem{corollaries}[theorem]{Corollaries}
\newcommand{\bth}{\begin{theorem}}
\newcommand{\bpr}{\begin{proposition}}
\newcommand{\epr}{\end{proposition}}
\newcommand{\bco}{\begin{corollary}}
\newcommand{\eco}{\end{corollary}}
\newcommand{\ble}{\begin{lemma}}
\newcommand{\ele}{\end{lemma}}
\newcommand{\bde}{\begin{definition}\rm}
\newcommand{\ede}{\end{definition}\rm}
\newcommand{\bre}{\begin{remark}\rm}
\newcommand{\ere}{\end{remark}}
\newcommand{\bex}{\begin{example}\rm}
\newcommand{\eex}{\end{example}}
\newcommand{\bcors}{\begin{corollaries}\rm}
\newcommand{\ecors}{\end{corollaries}}

\def\la#1{\hbox to #1pc{\leftarrowfill}}
\def\ra#1{\hbox to #1pc{\rightarrowfill}}
\def\fract#1#2{\raise3pt\hbox{$ #1 \atop #2 $}}

\def\lrar{{\ra 2}}

\def\sp#1{\hbox{SP}^{#1}}
\def\cp#1{\hbox{CP}^{#1}}
\def\gp#1{\hbox{$\Gamma$P}^{#1}}
\def\spy{\hbox{SP}^{\infty}}
\def\bsp#1{\overline{\hbox{SP}}^{#1}}
\def\sn{{\mathfrak S_n}}

\def\bbz{{\mathbb Z}}

\def\bbp{{\mathbb P}}
\def\bbr{{\mathbb R}}

\def\bbc{{\mathbb C}}
\def\bbr{{\mathbb R}}
\def\bbq{{\mathbb Q}}

\def\bbt{{\mathbb T}}

\begin{document}

\title{The Geometry and Fundamental Group
  of Permutation Products and Fat Diagonals}
  \author{Sadok Kallel and
  Walid Taamallah} \address{First author: Laboratoire Painlev\'e, Universit\'e des
  Sciences et Technologies de Lille, France,
  and American University of Sharjah, UAE.}\email{sadok.kallel@math.univ-lille1.fr}
  \address{Second author: Facult\'e des Sciences de Tunis}

\maketitle

\begin{abstract} Permutation products and their various ``fat diagonal"
  subspaces are studied from the topological and geometric
  point of view. We describe in detail the stabilizer and orbit stratifications
  related to the permutation action, producing a sharp upper bound for its depth
  and then paying particular attention to the geometry of the diagonal stratum.
  We write down an expression for the fundamental group
  of any permutation product of a connected space $X$ having the homotopy
  type of a CW complex in terms of $\pi_1(X)$ and $H_1(X;\bbz )$. We
  then prove that the fundamental group of the configuration space of
  $n$-points on $X$, of which multiplicities do not exceed $n/2$, coincides
  with $H_1(X;\bbz )$.
  Further results consist in giving conditions for when fat diagonal subspaces of
  manifolds can be manifolds again. Various examples and homological calculations
  are included.
\end{abstract}


\section{Introduction}

In this paper, the first in a sequel, we answer several basic questions about
the topology and geometry of permutation products and various analogs of
configuration spaces. These are spaces of relevance to all of algebraic and
geometric topology.

Let $\Gamma$ be a subgroup of the $n$-th symmetric group ${\mathfrak S}_n$,
and define the \textit{permutation product} $\gp{n}(X)$ to be the quotient of
$X^n$ by the permutation action of $\Gamma$ on coordinates.
The prototypical example being of course the $n$-th
symmetric product $\sp{n}(X)$ which corresponds to $\Gamma = {\mathfrak
  S}_n$. The cyclic product $\cp{n}(X)$ corresponds on the other hand to
$\Gamma = \bbz_n$ the cyclic group.  These spaces become interesting very
quickly as is illustrated for example with $\cp{3}(T):= T^3/\bbz_3$; the
cyclic product of order $3$ of an elliptic curve $T$, which has the structure
of an explicit bundle over $T$ (see Proposition \ref{elliptic}).

Permutation products have of course been studied early on, and a functorial
description of their homology described in \cite{dold}. The topology of the
spaces $\gp{n}(X)$ depends as one might expect on the way the group $\Gamma$
embeds in $\sn$. These are stratified
spaces \cite{hughes}. In general, the action of a finite group $\Gamma$ on a
space $Y$ stratifies both $Y$ and its quotient $Y/\Gamma$ according to
\textit{stabilizer} and \textit{orbit} types respectively, and
the bulk of the present work is understanding these stratifications
in the case of $\Gamma$ acting on $X^n$ by permutations.

In \S\ref{strats} we discuss the notion of ``depth" of a group
$\Gamma\subset\sn$ acting by permutation on $X^n$ and then compare it to the
length $\ell (\Gamma)$ of the group which we recall is
the length of the longest chain of subgroups in $\Gamma$.
The depth is bounded above by $\ell (\Gamma)$ and the
interesting part is to see that this upper-bound is attained for permutation
products derived from the permutation representation of $\Gamma$ (Theorem
\ref{depthgroup}).
In the case of a permutation product of a manifold, we can refine our
understanding of the associated strata.  We give in \S\ref{diagonalstratum} a
handy description of the tubular neighborhood of the diagonal $X$ in $\gp{n}X$
as a fiberwise cone on a fiberwise quotient bundle associated to a sphere
bundle over $X$ (Proposition \ref{neighdiagonal}).

We next turn to computing the fundamental group. We will assume throughout that
$X$ is of the homotopy type of a CW complex.  We recall that a permutation
subgroup $\Gamma\subset\mathfrak S_n : =\hbox{Aut}\{1,2\,\ldots,
n\}$ is called transitive if any $i$ is mapped to any $j$ via the action.
The following prototypical result is derived
in \S\ref{fundgroup}: \textit{If
  $\Gamma$ is a transitive subgroup of $\mathfrak S_n$ acting on $X^n$ by
  permutation of coordinates, and $n\geq 2$, then $\pi_1(\gp{n}X)$ is abelian and there is an
  isomorphism $\pi_1(\gp{n}X)\cong H_1(X;\bbz )$}. Very early on, P.A. Smith
in \cite{smith} has shown that $\pi_1(\sp{2}X)$ is abelian and his simple
geometric argument extends in an elementary and self-contained
manner to cover the more general case (see \S\ref{fundgroup}). Moreover since
any non transitive subgroup of $\sn$ is up to conjugation contained in a
standard Young subgroup ${\mathfrak S_k}\times {\mathfrak S_{n-k}}$,
permutation products of non transitive subgroups can be reduced to transitive
ones of smaller size and a general statement for $\pi_1(\gp{n}X)$ can be
given as follows\footnote{This result is motivated by a question of Mathai Varghese to the first author.}. Let's write $\Omega:=\{1,\ldots, n\}$ and again
$\sn = \hbox{Aut}(\Omega )$. We say an orbit is ``non-trivial" if it is not reduced to a single point.

\bth\label{main0} If $\Gamma\subset\sn$, then
$$\pi_1(\gp{n}X)\cong \pi_1(X)^{n-\sum_i k_i}\times H_1(X;\bbz )^r$$
where $r$ is the number of non-trivial transitive orbits under the action of
$\Gamma$ on $\Omega$, and $k_1,k_2,\ldots, k_r$ are the sizes of these
non-trivial orbits.
\end{theorem}

In the case $\Gamma=\mathfrak S_n$, we can further compare the higher
homotopy groups $\pi_i(\sp{n}X)$ to the corresponding homology
groups, and these coincide provided $0\leq i\leq r+2n-1$, where $r\geq 1$ is the
connectivity of $X$ and $n>1$ (Theorem \ref{stable}). This recovers an original
result of Dold and Puppe \cite{dp}.

The second part of this paper deals with the fat diagonal subspaces and their complements. For a given positive integer $d$, define $B_d(X,n)$
to be the subspace of $\sp{n}X$ of all unordered tuples
$[x_1,\ldots, x_n]$, or \textit{configurations}, such that at least $d$ entries
are equal.  These are the spaces we refer to as ``fat diagonals".
We have a decreasing
filtration
\begin{equation}\label{decfil} B_1(X,n) = \sp{n}(X)\supset
  B_2(X,n)\supset\cdots\supset B_n(X,n)=X \end{equation}
which interpolates between the symmetric product and the thin diagonal.
For example $B_2(X,n)$ is the standard fat diagonal in $\sp{n}X$
consisting of configurations with at least two points coinciding. In \S\ref{configs}
we characterize $B_d(-,n)$ as homotopy functors with certain properties, and then
give precise conditions for $B_d(X,n)$ to be a manifold if $X$ is
(Theorem \ref{manifold}).

As in the case of permutation products, it turns out that the fundamental
group of fat diagonals abelianizes ``whenever it can". The following is proved in
\S\ref{fundgroupbd} using an adaptation of Smith's argument.

\bth\label{main1} Let $X$ be a based connected topological space.
Then $\pi_1(B_d(X,n))\cong H_1(X;\bbz )$
provided that $1\leq d\leq {n\over 2}$.
\end{theorem}

Note that the case $d > n/2$ presents no particular interest since in this
case $B_d(X,n)\cong X\times\sp{n-d}X$.  More generally one shows that if
$\Gamma$ is a ``$d$-transitive" subgroup of $\sn$,
then $\pi_1(F_d (X,n)/\Gamma)$
is abelian for $n\geq 2$ and $1\leq d\leq {n\over 2}$ (Corollary \ref{dtransitive}).
Here $F_d(X,n)\subset
X^n$ is the so-called ``$d$-equal subspace" consisting of all tuples
$(x_1,\ldots, x_n)$ such that $x_{i_1}=\cdots = x_{i_d}$ for some choice of
sequence $i_1<\cdots < i_d$.

The Euler characteristics of the spaces $B_d(X,n)$ and their complements
$B^d(X,n):=\sp{n}(X) - B_{d+1}(X,n)$ will appear in upcoming work.
The fundamental group of $B^d(X,n)$ is determined in upcoming work by the first author and In\`es Saihi.

\vskip 5pt
\noindent{\sc Remerciements}: Nous tenons \`a remercier le PIms \`a Vancouver
(en particulier A. Adem et I. Ekeland), le CNRS (en particulier
J.M. Gambaudo) et le groupe de topologie de la facult\'e des sciences de
Tunis (en particulier Sa\"{\i}d Zarati et Dorra Bourguiba), pour avoir rendu
ce travail possible. Nous remercions \'egalement Dimitri Markushevich et
Zinovy Reichstein pour l'aide qu'ils nous ont apport\'ee dans la
d\'emonstration de la Proposition \ref{elliptic} et dans la discussion de la
\S\ref{strats} respectivement. Finalement nous remercions Bruce Hughes pour
nous avoir indiqu\'e la r\'ef\'erence \cite{beshears}.


\section{A geometric sample}\label{sample}

All spaces in this paper are assumed to be
first countable, locally compact, Hausdorff and path-connected.
Elements of $\gp{n}(X)= X^n/\Gamma$, with $\Gamma\hookrightarrow\sn$, are
written as equivalence classes $[x_1,\ldots, x_n] := \pi (x_1,\ldots, x_n)$,
where $\pi : X^n\lrar\gp{n}(X)$ is the quotient map. It will be convenient to call an
element in $\gp{n}(X)$ a \textit{configuration}. When $n=2$, the only
permutation product of interest is the symmetric square $\sp{2}(X)$.

The homotopy type of a permutation product $\gp{n}X$ depends on the way
$\Gamma$ embeds in $\sn$. Take for example $\Gamma = \bbz_2$ and embed it in
$\mathfrak S_4$ by sending the generator to either $(12)$ or to
$(13)(24)$. Then the permutation products corresponding to these embeddings
are respectively
$$\Gamma_1\hbox{P}^4(X) = \sp{2}(X)\times X\times X\ \ \
\hbox{and}\ \ \ \Gamma_2\hbox{P}^4(X) = \sp{2}(X\times X)$$
Of course conjugate embeddings yield homeomorphic quotients.

\bex\label{morton} Permutation products of the circle are particularly interesting.
They are special cases of what are called
``toroidal orbifolds" (see \cite{adem}).
Following an argument of Morton [22], we can see that for $\Gamma\subset\mathfrak S_n$, there is a bundle projection $\gp{n}(S^1)\lrar S^1$ with
fiber identified with the quotient $V/G_0$ where $V$ and $G_0$ are defined
as follows. First let $\Gamma$ and $\bbz^n$ act on $\bbr^n$ by permutations and
translations respectively. Define $G$ to be the affine subgroup generated by
 $\Gamma$ and $\bbz^n$. Then $V$ is the hyperplane
$\{x_1+\cdots +x_n=0\}\subset\bbr^n$ and $G_0\subset G$ is the stabilizer of $V$;
$G/G_0\cong\bbz$.
In the case of symmetric products, Morton shows that the quotient $V/G_0$
is an $(n-1)$-simplex, and so the product map $\sp{n}(S^1)\lrar S^1$ is a
disk bundle which is furthermore orientable if and only if $n$ is odd \cite{morton}. For
example $\sp{2}(S^1)$ is the M\"obius band, and the diagonal copy of $S^1\subset\sp{2}(S^1)$
is precisely the boundary of the band. In the case of the cyclic product
$\cp{3}(S^1)$ for example, the fiber of $\cp{3}(S^1)\lrar S^1$ is a copy of $S^2$
(Proposition \ref{c3} gives a complete treatment).
\eex

\bex\label{mattuck} Let $T=S^1\times
S^1$ be a torus, $S^1\subset\bbc$. Here too we have an (abelian) multiplication
$m: \sp{2}T\lrar T$ which is a bundle
projection and the fiber at $1$ is the quotient $T/\mathfrak S_2$, where
$\mathfrak S_2$ acts by taking $(x,y)\mapsto (\bar x,\bar y)\in T$. This action
corresponds to the hyperelliptic involution on the torus with $4$ fixed
points, and the quotient is a copy of the Riemann sphere.  This is of course
a special case of the well-known fact that multiplication $\mu: \sp{n}(T)\lrar T$ is
an analytic fiber bundle over the ``Jacobian" $T$ with fiber $\mu^{-1}(1)\cong\bbp^{n-1}$. \eex

When $n\geq 3$, the cyclic and symmetric products start to differ.  We can for
example use Morton's description of $\sp{n}(S^1)$ discussed in Example \ref{morton} to deduce
the following nice characterization \cite{wagner}.

\bpr\label{c3} There is a homeomorphism $\cp{3}(S^1) = S^1\times
S^2$. \epr

\begin{proof} The key point is the existence of two distinct sections to the
  projection $\pi: \cp{3}(S^1)\lrar\sp{3}(S^1)$. Indeed the clockwise (or
  counterclockwise) order of any three points on the circle is not changed by
  a cyclic $\bbz_3$-permutation which means that given $\{x_1,x_2,x_3\}\subset
  S^1$, any choice of writing these points in clockwise
  (resp. counterclockwise) orientation gives a well-defined element in
  $\cp{3}(S^1)$. The images of these two sections give two copies of
  $\sp{3}(S^1)$ in $\cp{3}(S^1)$. Since any triple of points in $S^1$ has
  either a clockwise or counterclockwise orientation, we have the
  decomposition
$$\cp{3}(S^1) = \sp{3}(S^1) \cup \sp{3}(S^1)$$
and the union is over the fat diagonal $B_2(S^1,3)$ of points of the form
$[x,x,y]$.  But $\sp{3}(S^1)$ is the product bundle over $S^1$ with fiber the
two-disk $D^2$ according to Morton's description (Example \ref{morton}),
and the configurations with repeated points consist of its
sphere bundle $\partial D^2\times S^1=S^1\times S^1$. This means that
$\cp{3}(S^1)$ is a union of a pair of disk bundles joined along their common
circle bundle; i.e. this is the \textit{fiberwise} pushout of two trivial disk bundles
$$S^2=D^2\cup_{S^1}D^2\lrar\cp{3}(S^1)\lrar S^1$$
and it is trivial as a bundle.
\end{proof}

It is slightly harder to determine $\cp{3}(T)$ when $T=S^1\times S^1$ is the
torus. There is a branched degree two covering $\pi : \cp{3}(T)\lrar\sp{3}T$
but there is no obvious section this time.  The following characterizes
$\cp{3}(T)$ algebraically and is to be compared to (\cite{bl}, Corollary 3.4).
We write the group structure on $T$ additively.

\bpr\label{elliptic} Let $T$ be an elliptic curve embedded in $\bbp^2$.
Addition $\mu_c : \cp{3}(T)\lrar T$,
$[x_1,x_2,x_3]\mapsto x_1+x_2+x_3$, is a
bundle projection with fiber a simply-connected
algebraic surface ramified over $\bbp^2$ along a sextic curve. The cohomology groups
of $\cp{3}(T)$ are given according to
$$
\begin{tabular}{|c|c|c|c|c|c|c|c|c|}
  \hline
  $H^{0}$&$H^{1}$& $H^{2}$ & $H^{3}$&$H^{4}$& $H^{5}$ & $H^{6}$\\
  $\bbz$&$\bbz^2$&$\bbz^5$&$\bbz^8$&$\bbz^5$&$\bbz^2\oplus\bbz_3$&$\bbz$ \\ \hline
\end{tabular}$$\epr

\begin{proof}
 The map $\mu_c$ makes $\cp{3}T$ into a bundle over $T$.
 The preimage $\mu_c^{-1}(0)$ is all equivalence classes of
 triple of points $[x,y,z]\in\cp{3}(T)$ with $x+y+z=0$.
 This can be easily identified with the quotient $\mu_c^{-1}(0) = (T\times T)/\bbz_3$ where
  $\bbz_3$ acts via its generator $\tau$ of order three as follows
$$\tau : (x,y)\mapsto (-x-y, x)$$
For simplicity we write $\bbt := T\times T$. We then have the map of bundles
$$\xymatrix{
  \bbt/\bbz_3\ar[r]\ar[d]^{2:1}&\cp{3}(T)\ar[r]^{\ \ \mu_c}\ar[d]^\pi\ar[d]&T\ar[d]^=\\
  \bbt/\mathfrak S_3\ar[r]&\sp{3}(T)\ar[r]^{\ \ \mu_s}&T }
$$
The symmetric group $\mathfrak S_3$ acts on $\bbt$ via its
generators $\tau$ of order $3$ and $J : (x,y)\mapsto (y,x)$ of order $2$. Of
course $\mu_c^{-1}(0)\rightarrow\mu_s^{-1}(0)$ is a branched degree two cover
(denoted by the usual notation $2:1$).
The branching locus is the image of the diagonal $C := \{[x,x], x\in T\}$ in
$\bbt/\mathfrak S_3$. This locus can now be explicitly identified. To that end
we start by identifying the quotient
 $\mu_s^{-1}(0)=\bbt/\mathfrak S_3$ with $\bbp^2$.
This we already know from Remark \ref{mattuck} but we need an explicit correspondence.
Let $(\bbp^2)^*$ be the dual of $\bbp^2$ consisting of projective lines in $\bbp^2$,
 and this is again a copy of $\bbp^2$. Since $T\subset\bbp^2$, consider the map
\begin{eqnarray*}
\bbt = T\times T\lrar (\bbp^2)^*\ \ ,\ \
(x,y)&\mapsto& \hbox{unique line in $\bbp^2$ through $x$ and $y$ if $x\neq y$}\\
(x,x)&\mapsto&\hbox{tangent to $T$ at $x$}
\end{eqnarray*}
This map is surjective since any line in $\bbp^2$ must intersect $T$,
and it descends to the quotient $\bbt/\mathfrak S_3$ simply
because by the very definition of the group law on $T$, the line through $(x,y)$
meets again the curve at the point $-x-y$ so that $(x,y)$ and $(-x-y,x)$ define
the same line in $\bbp^2$. This leads to the identification
$\mu_s^{-1}(0)\cong (\bbp^2)^*$.
The branch locus curve $C$ maps isomorphically to
the dual curve $T^*$ made out of the tangent lines to $T$.  But
the dual curve to a smooth plane cubic has 9 singular points (cusps) which
correspond under the duality to the 9 flexes of $T$ (i.e the inflection points
of $T$ or also the points of order $3$ in the group law for $T$).  Note that this dual
curve is of course of genus $1$ and hence must have degree $6$ (i.e. a sextic)
by the degree-genus formula.

To summarize, the map $\mu_c^{-1}(0)\lrar\mu_s^{-1}(0)$ is identified with a
map $\mu_c^{-1}(0)\lrar\bbp^2$ which ramifies over the curve $T^*$.
Over the cusps in $T^*$ lie the  $9$ singular points of $\mu_c^{-1}(0)$.
Note that $\mu_c^{-1}(0)$ is simply connected because of the short exact sequence
$$ 0=\pi_2(T)\lrar\pi_1(\mu_c^{-1}(0))\lrar\pi_1(\cp{3}T)\fract{\mu_{c*}}{\lrar}\pi_1(T)\lrar 0$$
associated to the bundle $\mu_c :\cp{3}T\lrar T$ and the fact that $\mu_{c*}:
\pi_1(\cp{3}T)\lrar\pi_1(T)$ must be an isomorphism between two copies of
$\bbz\oplus\bbz$ according to Theorem \ref{main0} next and the existence of a section.

Finally and for any prime $p$, the cohomology of $\cp{p}(T)$ can be deduced from \cite{adem}.
Here one thinks of $\cp{3}(T)$ as a ``toroidal orbifold", that is
as a quotient of $(\bbr^2/\bbz^2)^3$ by an
action of $\bbz_3$ preserving the lattice. This action is the sum of two
regular representations $\bbz_3\lrar GL_3(\bbr )$ since any such
representation corresponds to cyclically permuting the basis vectors.
The calculation above for $H^*(\cp{3}(S^1\times S^1))$ is then deduced from
(\cite{adem}, Theorem 1). Details in \cite{taamallah}.
\end{proof}


\section{Stabilizer and Orbit stratifications}\label{strats}

The category of $G$-stratified spaces is a suitable category in which to
consider permutation products.  As in \cite{beshears}, we assume $G$ to act
 discontinuously\footnote{\cite{brown}, Definition 1.3:
A group $G$ acts discontinuously on a
space $X$ if the stabilizer of each point of $X$ is finite, and each point
$x\in X$ has a neighborhood $V_x$ such that any element $g$ of $G$ not in the
stabilizer of $x$ satisfies $V_x\cap g\cdot V_x = \emptyset$.} and we define
for a $G$-space $X$, a stratum for each conjugacy
class $(H)$ of a subgroup $H\subset G$ as follows
$$X_{(H)} := \{x\in X\ |\ G_x\sim H\} = \bigcup_{K\sim H}X_K$$
where $\sim$ means ``being conjugate to'', $G_x:=\{g\in G\ |\ gx = x\}$
the isotropy group of $x$, and
$$X_K := \{x\in X\ |\ G_x = K\}$$
In words, $X_K$ consists of those points $x\in X$ that are fixed by $K$ and by
no larger subgroup of $G$. This space is contained in the fixedpoint set of $K$
and it depends strongly on $G$ as well.
We can see that $X_{(G)}=X_G$ consists of the fixed points of the action.
On the other hand $X_{(Id)}=X_{Id}$ is
the subset of points of $X$ on which $G$ acts freely; here $Id$ is the trivial subgroup.
As $H$ ranges
over all subgroups of $G$ we get a stratification of $X$ (i.e. a partition of
$X$, see \cite{hughes}) called the
\textit{stabilizer stratification} \cite{beshears}(\footnote{We depart from the
terminology of Beshears who calls this stratification the ``orbit type" stratification.}).
The \textit{strata} are the various $X_{(H)}$, $H\subset G$.
One checks that $X_{(H)}$ are
$G$-invariant subspaces since $G_{gx} = g\cdot\ G_x\cdot g^{-1}$, and that the
corresponding quotients $X_{(H)}/G$ stratify $X/G$ as well making up the so-called
\textit{orbit stratification}.

\bex The stabilizer stratification of $\bbz_2$ acting on $S^n$ by reflection
with respect to the plane of an equator has two strata : the two open
hemispheres form a stratum (the ``free part") and the equator forms another
(the ``fixed point" part). \eex

\bex Strata corresponding to proper subgroups can be empty.  For example if
$G$ is a non-trivial group acting trivially on a space $X$; i.e. $gx=x$ for
all $g\in G$, $x\in X$, then $X_{(H)} = \emptyset$ for all proper subgroups
$H\subset G$. In this case there is only one (non-empty) stratum
consisting of the fixed points $X_{(G)}=X$.\eex

\noindent {\sc General Properties}: Assume $G$ is acting on $X$, $H\subset G$ a
subgroup. We will write throughout $\overline{A}\subset X$ the closure of $A$ in $X$.
\begin{itemize}
\item If we set $H^g=g^{-1}Hg$, then $X_{H^g}=g^{-1}X_H$ and so we can write
$X_{(H)} = \bigcup_{g\in G}gX_H$.
\item If $x\in {\overline X}_H$, then $H\subset G_x$ the stabilizer of $x$ in $G$.
  Being in the closure of $X_{H}$, there is a sequence of
  points $x_i\in X_{H}$ converging to $x$ (all spaces are first countable as mentioned). But for any $h\in H$, $hx_i=x_i$ and so $hx_i$ converges to $hx = x$, and $H\subset G_x$.
\item If $X_{K}\cap {\overline X}_{H}\neq\emptyset$, then $H\subset K$. This is
a direct consequence of the preceding property. The
  reciprocal statement is not always true: let $\mathfrak S_4$ act on $X^4$ by
  permuting coordinates, $H=\bbz_4\subset \mathfrak S_4=K$. Then
  $X^4_{(\bbz_4)}$ is empty but $X^4_{(\mathfrak S_4)}\cong X$ is not.
\end{itemize}

For a stratified space $X$ one defines the relation $\leq $ on the index set
$\mathcal I$ of a stratification by setting $i\leq j$ if and only if
$X_i\subset \overline{X}_j$; the closure of $X_j$ in $X$.  A stratification
$\{X_i\}_{i\in\mathcal I}$ satisfies the ``frontier condition" if for every $i,j\in\mathcal I$
$$X_i\cap \overline{X}_j\neq
\emptyset\ \ \ \hbox{implies}\ \ \ \ X_i\subseteq \overline{X}_j$$
(and $i<j$ if $i\neq j$).
For example a locally finite simplicial complex is such a
stratified space with strata the open simplices. Similarly if $N$ is a
submanifold of $M$, then $N$ and $M-N$ form a stratification of $M$ satisfying
the frontier condition.

It is not hard to come up with group actions where both
the orbit and stabilizer stratifications do not satisfy the
frontier condition. For example $X$ is a fork with $3$ tines and $G=\bbz_2$ acting
by spinning the fork around its axis an angle $\pi$. The following lemma gives a
sufficient condition for when such stratifications have the frontier condition.

\ble\label{sufficient} Let $\Gamma$ be a finite group acting on $X$, and suppose
that for any pair of subgroups $H\subset K$ of $\Gamma$ such that
$X_H\neq\emptyset$, we have $X_K\subset \overline{X}_H$.
Then the stabilizer stratification satisfies the frontier condition. Moreover
\begin{equation}\label{union}
\overline{X}_{(H)} = \bigcup_{K\supset
  H}X_{(K)}
  \end{equation}
\ele

\begin{proof}
  Suppose that $X_{(K)}\cap \overline{X}_{(H)}\neq\emptyset$. Then
  $X_K\cap\overline{X}_{H^\sigma}\neq\emptyset$ for some conjugate subgroup
  $H^\sigma$ of $H$.  By one of the general properties mentioned above,
  $H^\sigma\subset K$. Since $X_H\neq\emptyset$, then
  $X_{H^\sigma} = \sigma^{-1}X_H$ is also non-empty and hence by our
  hypothesis $X_K\subset\overline{X}_{H^\sigma}$ so that
  $X_K\subset\overline{X}_{(H)}$ and then obviously
  $X_{(K)}\subset\overline{X}_{(H)}$.
  This shows that the frontier condition is satisfied.
  To get the claimed description of $\overline{X}_{(H)}$, one needs to see first that
  $\bigcup_{K\supset H}X_K\subset\overline{X}_H$ and this is an immediate consequence of the
  hypothesis. The reverse inclusion
  $\overline{X}_H\subset \bigcup_{K\supset H}X_K$
  follows from one of properties listed earlier
  that if $x\in \overline{X}_H$, then $H\subset K$ where $K$ is the stabilizer of $x$.
  We then have that $\overline{X}_H= \bigcup_{K\supset H}X_K$.
  Since on the other hand ${X}_{(H)} = \bigcup_{H'\sim H}{X}_{H'}$, we can pass to
  the closure and the equality (\ref{union}) follows.
\end{proof}

\subsection{Permutation Products}\label{structure}

In the case of permutation products, the orbit and stabilizer stratification
turn out to be well behaved. The first point to make is that
for path-connected $X$,
the stabilizer of a point $(x_1,\ldots,x_n)\in X^n$ only depends on the multiplicities
of the $x_i$'s and not on their location in $X$. Let $\sigma\in\Gamma\subset\sn$ be
a permutation. It has a
\textit{cycle structure}; i.e. can be expressed as a product of pairwise
disjoint cycles
\begin{equation}\label{product}
\sigma = (i^1_1\ldots i^1_{d_1})(i^2_{1}\ldots i^2_{d_2})\cdots (i^k_{1}\ldots i^k_{d_k})
\end{equation}
Fixed points of $\sigma\in\Gamma$ with such a cycle structure
consist of all $n$-tuples having equal entries at $i_1^1,...,i^1_{d_1}$, then
equal entries at $i_1^2,\ldots, i^2_{d_2}$, etc.

Consider the case of $\sn$ acting on $X^n$. Throughout we will think of $\sn = \hbox{Aut}(\Omega)$ as the permutation group of $\Omega:=\{1,\ldots, n\}$.
Let $\mathfrak S_{d_1}\times\ldots \times
\mathfrak S_{d_k}$ be the subgroup where $\mathfrak S_{d_1}=\hbox{Aut}\{1,\ldots, d_1\}$,
$\mathfrak S_{d_2}=\hbox{Aut}\{d_1+1,\ldots, d_1+d_2\}$, etc, with $d_i\geq 1$ and $\sum d_i=n$.
We refer to \textit{Young} any
subgroup $H\subset\sn$ conjugate to such a product. The
stratum $X^n_{(\mathfrak S_{d_1}\times\ldots \times \mathfrak S_{d_k})}$
corresponds to all $n$-tuples of points which after permutation can be brought
to the form
\begin{equation}\label{fixedpoint}
  (\underbrace{x_1,\ldots, x_1}_{d_1}, \underbrace{x_2,\ldots, x_2}_{d_2}, \cdots, \underbrace{x_k,\ldots, x_k}_{d_k})
\end{equation}
with $x_i\neq x_j$ for $i\neq j$ and $\sum d_i = n$.   From this
description, it is easy to see that $X^n_K$ is trivial if $K$ is not a Young
subgroup of $\sn$.  Moreover each cycle structure (\ref{product}) corresponds to a unique conjugacy class of a Young subgroup of $\sn$. This implies in particular that for
various choices of Young subgroups, the subspaces $X^n_{(young)}$ give rise to
a stratification of $X^n$, which is moreover $\mathfrak S_n$-invariant. We refer to
the subspace $X^n_{(young)}$ as a ``Young stratum".  The closure of the Young
stratum $X^n_{(\mathfrak S_{d_1}\times\ldots \times \mathfrak S_{d_k})}$ is
made up of all Young strata obtained by merging $\mathfrak S_{d_1}\times\ldots
\times \mathfrak S_{d_k}$ into Young subgroups with smaller factors, so for
example $X^n_{(\mathfrak S_{d_1+d_2}\times\mathfrak S_{d_3}\ldots \times
  \mathfrak S_{d_k})}$ is in the closure. In fact this closure corresponds to the fixed point sets since
  $$\overline{X^n}_{\mathfrak S_{d_1}\times\ldots
\times \mathfrak S_{d_k}} = Fix(\mathfrak S_{d_1}\times\ldots
\times \mathfrak S_{d_k})$$
To sum up, $X^n$ under the
  action of $\sn$ is stratified by the Young strata and when $X$ is path
  connected, this stratification satisfies the frontier condition.
The latter fact turns out to be true for any finite group acting by permutation
on a path-connected space as we now explain.
Before we do so, we need the following useful remark.

\bre\label{fixedset} We give here a complete description of
Fix$(H)\subset X^n$ for any subgroup $H\subset\sn$.
View again $\sn =$ Aut$(\Omega )$,
$\Omega = \{1,2,\ldots, n\}$. The action of $H$ on $\Omega$ breaks it into orbits
$\Omega = \Omega_1\sqcup\cdots\sqcup\Omega_d$ of sizes $n_1,\ldots, n_d$ respectively,
with $\sum n_i=n$. Then $(x_1,\ldots, x_n)\in X^n$
is fixed by $H$ if and only if $x_i=x_j$ whenever
$i,j\in\Omega_r$ for some $r$. To each orbit $\Omega_j$ corresponds a subtuple
of $(x_1,\ldots, x_n)$, or a ``bloc", having equal entries.
For example if $H=\bbz_2\times\bbz_3\subset\mathfrak S_5$, $\bbz_2$ generated by $(12)$ and
$\bbz_3$ generated by $(345)$, then
$$
\hbox{Fix}(H) = \{(x,x,y,y,y)\in X^5\}
$$
The orbits are
$\Omega_1=\{1,2\}$ and $\Omega_2=\{3,4,5\}$. To $\Omega_1$ corresponds the ``bloc"
of $x$'s and to $\Omega_2$ the bloc of $y$'s.
\ere

\bpr\label{frontier} Suppose $X$ is path-connected and let
$\Gamma\subset\sn$ act on $X^n$ by permutation. Then both the associated
stabilizer and orbit stratifications satisfy the frontier condition.
\epr

\begin{proof}
Let $H\subset\sn$ and let
${\mathcal P}_H$ be the partition of orbits $\Omega_1,\ldots, \Omega_d$
associated to $H$ as in Remark \ref{fixedset}.
Then a tuple $(x_1,\ldots, x_n)\in X^n$ is fixed by $H$ if it decomposes into blocs
of equal entries, each bloc determined by some $\Omega_j$.
We will say in this case that $(x_1,\ldots, x_n)$
is ``subordinate" to the partition ${\mathcal P}_H$.
Notice that a tuple $(x_1,\ldots, x_n)$ is in the closure $\overline{X^n}_H$ if
$x_i=x_j$ whenever $i$ and $j$ are in some union of the $\Omega_r$'s; that is whenever it is
subordinate to a coarser partition.

Assume that $H\subset K$ two subgroups of $\Gamma$. Then the
partition of orbits of
$H$ is finer then that defined by $K$. For path-connected $X$, it is easy to
arrange that any $(x_1,\ldots, x_n)$ subordinate to a partition
${\mathcal P}$ is the limit of tuples subordinate to a finer partition.
Assume $\overline{X^n}_H\neq\emptyset$.
Then any point in $X^n_K$ (if non-empty) must be a limit point in $X^n_H$ and so $X^n_K\subset \overline{X^n}_H$. According to Lemma \ref{sufficient} the stabilizer stratification satisfies the frontier condition.
\end{proof}

\subsection{Depth} We assume here that all groups are finite.  Given a finite
stratification $\{X_i\}_{i\in\mathcal I}$ of a space $X$ satisfying the
frontier condition, we define the \textit{depth} of a stratum $X_s$ to be the
maximal length $k$ of a sequence $s=s_0 < s_1 < \cdots < s_k$ in $\mathcal I$.
We recall that $s_i<s_j$ means $X_{s_i}\subset\overline{X}_{s_j}$.
The depth of a (finite) stratification as a whole is the maximum over the
depths of its strata.
In this section we assume that stabilizer stratifications associated to group
actions satisfy the frontier condition.

\bde Represent an action of $\Gamma$ on a space $X$ by a homomorphism $\phi
: \Gamma\lrar \hbox{Homeo}(X)$. Denote by $d_\phi (\Gamma ,X)$ the depth of the
stabilizer stratification associated to this action. As is clear, this depth
depends on both $\phi$ and $X$.  \ede

For given $\Gamma$, the various depths $d_\phi (\Gamma ,X)$ can be compared to
the length of $\Gamma$.  Recall that the length $\ell (\Gamma )$ of a finite
group $\Gamma$ is the length $\ell$ of a longest chain of subgroups
$\{1\}=\Gamma_0\subset \Gamma_1\subset\cdots \subset \Gamma_\ell= \Gamma$.

\bth\label{depthgroup} If $\Gamma$ is a finite group, then $0\leq
d_\phi(\Gamma,X)\leq \ell (\Gamma )$ and these equalities are sharp.
\end{theorem}

\begin{proof}
  To see that $d_\phi(\Gamma, X)$ can be $0$, simply choose $X$ to be a point.
  To show the right inequality, recall by one of our earlier properties,
  that in any stabilizer stratification,
  $\emptyset\neq X_{(K)}\subset {\overline X}_{(H)}$
  implies  that $H\subset K^\sigma$ for some conjugate subgroup $K^\sigma\sim K$.
  It follows that any chain of strata in the stratification
  must correspond to a chain of subgroups and hence the depth of the
  stratification is at most $\ell (\Gamma)$.
  To see that this inequality can be
  sharp, we will consider the permutation
  action of $\Gamma$ on the set of its elements $\Omega=\{1,\ldots, |\Gamma |\}$.
  This action is free and transitive. It
gives an embedding $\psi : \Gamma\hookrightarrow \mathfrak S_n$, $n=
  |\Gamma |$ (the so-called permutation representation of
  $\Gamma$) and an induced permutation action on $X^n$
$$\phi : \Gamma \lrar\ \hbox{Homeo}(X^n)\ \ \ ,\ \ \ n= |\Gamma |$$
We claim that $d_{\phi}(\Gamma, X^n)=\ell (\Gamma )$ for path-connected $X$.
The main step is to show that any subgroup $H\subset\Gamma$ gives rise to a
non-empty stratum $X^n_{(H)}$. We proceed by exhibiting an element in $\hbox{Fix}(H)$
which cannot be fixed by any larger subgroup. For ease we write $\psi (H)$ to refer
to $H$ acting on $\Omega$ and we write $H$ when we wish to indicate that
$H$ acts on $X^n$.

As before, $(x_1,\ldots, x_n)\in \hbox{Fix}(H)$ if it can be decomposed into $d$-blocs corresponding to $\Omega_1,\ldots, \Omega_d$; the orbits of the action of $\psi (H)$ on $\Omega=\{1,\ldots, n\}$. Each bloc is a subtuple characterized by the fact that its entries are equal. Moreover $H$ acts on each bloc transitively and doesn't permute elements from different blocs. We will choose $\zeta=(x_1,\ldots, x_n)\in \hbox{Fix}(H)$ such that
$x_i\neq x_j$ if they are not in the same bloc; i.e.
if $i,j$ not in the same $\Omega_r$ (of course this is possible if $X$ is not reduced to a point). We claim that $\Gamma_{\zeta}=H$. This is where the permutation
representation intervenes.

Let $g\in\Gamma_\zeta$.
Then $g$ leaves the bloc corresponding to $\Omega_1$ with entries equal to $x_1$ say
invariant. All orbits $\Omega_i$ have cardinality at least two since $\psi (H)$ acts
freely on $\Omega$. Choose $i,j\in\Omega_1$ such that $\psi (g)(i)=j$.
Since $H$ acts transitively on that bloc, there is $\psi (h)\in \psi (H)$ such
that $\psi (h)(i) =j$ as well.  This
means that at the level of the group $\Gamma$, there are generators
$g_i,g_j\in\Gamma$ such that $gg_i = g_j$ and $hg_i = g_j$, which means that
$h = g$ and $g\in H$. This shows therefore that $\zeta\in X^n_{H}$ and the stratum
$X^n_{(H)}$ is non empty.

The previous argument asserts that for every $H\subset \Gamma$, we have a non-empty
stratum $X^n_{(H)}$. Another argument of limit
using the path-connectedness of $X$ readily shows that if $H\subset K$, then $X^n_K\cap
\overline{X^n}_{H}\neq\emptyset$. It follows by Proposition \ref{frontier} that
$X^n_{(K)}\subset \overline{X^n}_{(H)}$.
Any chain of subgroups gives therefore a non-empty chain of strata
and $d_\phi (\Gamma, X^n)
\geq \ell (\Gamma )$, so that $d_\phi (\Gamma, X^n) = \ell (\Gamma )$ concluding
the proof.
\end{proof}

\bde For a given embedding $\phi :
\Gamma\hookrightarrow\sn$, we can write $d_\phi (\Gamma )$ the depth of the
orbit stratification of $\Gamma$ acting on $X^n$ where $X$ is path-connected and
not reduced to point.
This doesn't depend on $X$ as we pointed out. \ede

\bex  For $\phi: \Gamma\hookrightarrow\sn$, $d_{\phi}(\Gamma )\leq\ell (\Gamma)
 \leq\ell (\sn)$; the length of $\sn$
(Theorem \ref{depthgroup}). By a beautiful theorem of Cameron,
Solomon and Turull \cite{cst}, this is
$$\ell ({\mathfrak S}_n) = [(3n-1)/2] - b(n)$$ where $b(n)$ is the number
of ones in the dyadic representation of $n$.
\eex

\bex\label{symprod} Let $\phi=Id$ be the identity of $\sn$ with associated quotient space
$\sp{n}(X)$. The strata of the orbit type stratification are in one to
one correspondence with integer partitions of the form
$$P = [p_1,p_2,\ldots, p_k]\ \ , \ \ p_i\leq p_{i+1}\ ,\ \sum p_i =
n$$
To the partition
$[p_1,p_2,\ldots, p_k]$ corresponds the stratum of points of the form
$$[x_1,\ldots, x_1,x_2\ldots, x_2,\cdots ,x_k,\cdots x_k]$$
where $x_i\neq x_j$ if $i\neq j$, and each $x_i$ repeats $p_i$-times.
This is the stratum $X^n_{(\mathfrak S_{p_1}\times\ldots \times \mathfrak S_{p_k})}/\sn$.
The generic (open and dense) stratum is the configuration space of distinct points
$B(X,n) = F(X,n)/\sn = X^n_{(1)}/\sn$ and this corresponds to the partition
$[1,\ldots, 1]$. The diagonal in $\sp{n}X$ corresponds to the partition $[n]$.
The depth of this stratification is $d_{Id}(\sn) = n-1$.
The poset of strata here is the set of partitions partially ordered according to
$P\leq Q$ if $Q$ is a
refinement of $P$, that is if $Q$ is obtained by further partitioning the
parts of $P$.
\eex

\bex The depth of $\bbz_n$ acting on $X^n$ is the length of $\bbz_n$. This is because the
inclusion $\bbz_n\hookrightarrow\sn$, sending the generator to the $n$-cycle
$(12\cdots n)$ corresponds to the permutation representation of $\bbz_n$ (Theorem \ref{depthgroup}). Note
that if $p|n$, then $X^n_{(\bbz_p)}=X^n_{\bbz_p}$ in the stabilizer stratification associated to
 the action of $\bbz_n$ on $X^n$, consists up to cyclic permutation of all
points of the form
$$(x_1,x_2,\ldots, x_{n/p}, x_1,x_2,\ldots , x_{n/p},\ldots,
x_1,x_2,\ldots , x_{n/p})$$ with $x_i\neq x_j$ if $i\neq j$, each $x_i$
repeating $p$ times.  \eex


\section{Permutation Products of Manifolds}\label{diagonalstratum}

The fact that $\sp{n}(S)$ is a complex manifold, for $S$ a smooth algebraic
curve, is an exceptional property that is no longer true for general
$\gp{n}(S)$. For example here's a standard argument as to why $\cp{3}(\bbc )$ cannot
be homeomorphic to $\bbr^6$: decompose
$\bbc^3=L\times H$, $L$ the diagonal line in $\bbc^3$ and $H$ its orthogonal
hyperplane.  The linear action of $\bbz_3$ takes place entirely in $H$. This action
restricts to a fixed point free action on the unit sphere $S^3\subset H$
(since we have factored out the diagonal which consists of the fixed points).
It follows that $\bbr^6/\bbz_3\cong \bbr^2\times cM$ where $M=S^3/\bbz_3$ is a
three manifold with $H_1(M)=\bbz_3$ and $cM$ is the cone on $M$ (with apex at
the origin and the cone extending to $\infty$). This quotient cannot be
homeomorphic to $\bbr^6$.

In general, if $\Gamma$ is a finite group acting on a smooth manifold,
then its associated stabilizer
stratification is ``Whitney stratified'' (see \cite{pflaum}, Theorem 4.3.7).
In particular submanifold strata do have tubular neighborhoods of which boundaries
are not anymore sphere bundles.
We analyze here in details a neighborhood of the diagonal stratum
$M\hookrightarrow \gp{n}M$, $x\mapsto [x,\ldots, x]$ for $n\geq 2$.

Let $V$ be a (real) vector space
and let $\mathfrak S_n$ act on $V^n$ by permutation. The orthogonal complement
of the diagonal subspace is $H = \{(v_1,\ldots, v_n)\in V^n, \sum v_i = 0\}$ and
has by restriction an action of $\mathfrak S_n$.
Consider the map
$$\Phi : V^{n-1}\lrar V^n\ \ \ ,\ \ \
(v_1,\ldots, v_{n-1})\longmapsto (-v_1,\ldots, -v_{n-1}, v_1+\cdots
+v_{n-1})$$
If we require that this map be
$\mathfrak S_n$-equivariant, then the action on the left factor $V^{n-1}$ is determined
implicitly and we refer to it as the ``induced'' action.  For example for $n=2$
the induced action on $V$ is generated by the involution $v\longmapsto -v$.
If $n=3$, then the generating $3$-cycle $\sigma\in\mathfrak S_3$ acts on $V^3$ via
$\sigma (-v_1,-v_2,v_1+v_2) = (v_1+v_2,-v_1,-v_2)$ and this dictates
the induced action $\sigma (v_1,v_2) = (-v_1-v_2, v_1)$ in $V^2$
(same as in Proposition \ref{elliptic}).
The point is that the map $\Phi$ maps $V^{n-1}$ isomorphically and equivariantly onto $H$.

Let now $TM^{\oplus n-1}$ be the $(n-1)$-fold whitney sum of the tangent bundle
$TM$ of $M$, $\dim M=m$, and let $\Gamma\subset\mathfrak S_{n}$ act on each
fiber $V^{n-1}\cong (\bbr^{m})^{n-1}$, via the induced action.  The fiberwise
quotient is a bundle with fiber
$V^{n-1}/\Gamma$ homeomorphic to $H/\Gamma$. But $\Gamma$ acts
linearly on $H\subset V^n$ so that the quotient is the open cone $c(S/\Gamma)$,
where $S\cong S^{(n-1)m-1}$ is the unit sphere in $H$. The action being
fiberwise, we obtain therefore a new bundle $S_{\Gamma}M$ over $M$ with fiber
homeomorphic to $c(S/\Gamma)$.

\bpr\label{neighdiagonal} For $M$ a smooth closed manifold, a neighborhood
deformation retract $V$ of the diagonal $M$ in $\gp{n}M$ is homeomorphic to
the total space of $S_{\Gamma }M$.
The section corresponding to the cone points is the inclusion of $M$ into
$V$.  \epr

\begin{proof}
  The normal bundle of the diagonal copy of $M$ in $M^n$ is isomorphic to the
  Whitney sum $TM^{\oplus n-1}$. The reason for this is that if a point $x$ is
  sufficiently close to a point $y$, then it determines a unique vector $v\in
  T_xM$ such that $exp_x(v)=y$.
  So if $(x,y_1,\ldots, y_{n-1})$ is a point in $M^n$ with $y_i$ sufficiently
  close to $x$, then it determines
  a tuple of vectors $(v_1,\ldots, v_{n-1})$ with  $exp_x(v_i)=y_i$.
    For each $x\in M$, we write $T_xM^{\oplus
    n-1}$ the fiber at $x$ of $TM^{\oplus n-1}\lrar M$.

Choose a metric on
$M$ with injectivity radius greater than 2 and let $DM\subset TM$ be the
unit disc bundle of the tangent bundle of $M$. We will write $exp$ the
exponential map $DM\lrar M$ and $exp_x$ its restriction to the fiber $D_xM$
at $x \in M$.  Let $DTM^{\oplus n-1}$ be the unit disk bundle in the
Whitney sum. We get a map
$$DTM^{\oplus n-1}\fract{\phi}{\lrar} M^n$$
defined on each fiber by
$$(x, (v_1,\ldots, v_{n-1}))\longmapsto
(exp_x(-v_1), \cdots, exp_x(-v_{n-1}), exp_x(v_1+\cdots +v_{n-1}))
$$
This map is $\mathfrak S_n$-equivariant by the very definition of the
induced action, it is one-to-one and onto its image $V$.
Since $exp_x(0) = x$, this image produces a $\mathfrak S_n$-invariant
neighborhood of the diagonal in $M^n$ which is then homeomorphic to the total
space of $DTM^{\oplus n-1}$. A tubular
neighborhood of $M$ in $\gp{n}$ is therefore homeomorphic to the fiberwise
quotient of $DTM^{\oplus n-1}$ by the action of $\Gamma$ and this
is in turn homeomorphic to $S_{\Gamma}M$.
\end{proof}

\bco\label{neighborhood} Assume $M$ is a closed smooth manifold. A
neighborhood deformation retract $V$ of the diagonal $M$ in $\sp{2}M$ is
homeomorphic to the fiberwise cone on the projectivized tangent bundle of $M$.
\eco

\begin{proof} According to Proposition \ref{neighdiagonal},
the fiber at $x\in M$ of $S_{\bbz_2}M$ is the cone on $S^{m-1}/\bbz_2$
where $S^{m-1}$ is the unit tangent sphere at $x$ and $\bbz_2$ is acting
via the antipodal action $v\mapsto -v$.
\end{proof}

\bre In the case of the sphere, the complement of the diagonal in $\sp{2}S^n$
is $B(S^n,2)$ and this is known to deformation retracts onto $\bbr P^n$. \ere


\section{The Fundamental Group of Permutation Products}\label{fundgroup}

View $\sn$ as the permutation group of $\Omega=\{1,2,\ldots, n\}$. The set
$\Omega$ breaks into orbits under the action of $\Gamma$ viewed as a subgroup of $\sn$,
and we will say that an orbit is non-trivial if it contains at least two elements.
In this section we prove main Theorem \ref{main0} in the introduction. Let $r$ be the number of non-trivial transitive orbits under this permutation action
of $\Gamma$ on $\Omega$ and let $k_1,k_2,\ldots, k_r$ be the sizes of the \textit{non-trivial} orbits. If $\pi_1(X)$ is based at  $\ast\in X$, then $\pi_1(\gp{n}X)$ is based at $[\ast,\ldots,\ast]$.

\bth\label{mainstat} There is an isomorphism
$\pi_1(\gp{n}X)\cong \pi_1(X)^{n-\sum_i k_i}\times H_1(X;\bbz )^r$
\end{theorem}

We are aware that this theorem can be derived using the methods
of Brown and Higgins \cite{brown}. Our approach however is geometrically appealing
and totally self-contained. We start by proving that if $\Gamma$ is a transitive subgroup of $\mathfrak S_n$ acting by permutations on $X^n$ and $n\geq 2$, then $\pi_1(\gp{n}X) \cong H_1(X;\bbz )$(\footnote{We know of only one place where this result is mentioned (\cite{dp}, 12.15) and it is said there that this isomorphism in the case of a transitive action
would follow from earlier results of Dold, without further details.}).
Universal examples of transitive actions is when $H\subset \Gamma$ is a subgroup
and $\Gamma$ acts transitively on the right cosets of $H$ \cite{isaacs}.

\bex A special example of a transitive action
is when $\Gamma = \mathfrak S_\ell\wr\mathfrak S_k$ is a
wreath product (i.e.  these form the maximal, transitive and imprimitive
subgroups of $\mathfrak S_{\ell k}$, $\ell ,k>1$). In this case
$\gp{\ell k}(X)\cong\sp{k}(\sp{\ell }X)$ and
$\pi_1(\gp{\ell k}(X))\cong H_1(\sp{\ell }X;\bbz )\cong H_1(X;\bbz)$, consistently
with the fact that $\Gamma$ has only one transitive orbit.  \eex

In order to prove Theorem \ref{mainstat} we need some useful facts about
``path-lifting" of maps.
We say a map $\pi : Y\lrar X$ has the ``path-lifting property" if we can lift any path $[0,1]\lrar X$ to a path in $Y$ knowing where it starts in $Y$.
The following is Proposition 1.5 of \cite{brown}.

\bpr If $G$ is a group acting discontinuously
on a Hausdorff space $X$, then $\pi : X\lrar X/G$ has the path lifting property.
\epr

\bco\label{fixedpoint2} If a finite group $G$ acts
on a Hausdorff connected space $X$ with a fixed point $x_0$, then
$\pi_1(X,x_0)\lrar\pi_1(X/G,x_0)$ is surjective.  In particular
$\pi_1(X^n)\lrar\pi_1(\gp{n}X)$ is surjective for any $\Gamma\subset\mathfrak
S_n$.\eco

\begin{proof} A loop $\gamma$ based at $x_0$ in $X/G$ must lift to a loop
$\tilde\gamma$, and
the homotopy class $[\tilde\gamma]$ maps to $[\gamma ]$. This establishes the main statement. In the case of $X^n/\Gamma$,
one chooses as basepoint a diagonal element $[\ast,\ldots,\ast]$, $\ast\in X$.
\end{proof}

\bre Note that the surjection $\pi_1(X^n)\twoheadrightarrow\pi_1(\gp{n}X)$ is
equivalent to the following intuitive fact: if ${\bf x}=[x_1,\ldots,
x_n]$, ${\bf y} = [y_1,\ldots, y_n]$ are two configurations in $\gp{n}(X)$, and
$\gamma$ a path between them, then there are paths $\gamma_1,\ldots,
\gamma_n$ in $X$ such that
$\gamma (t) = [\gamma_1(t),\ldots, \gamma_n(t)]$ with $\gamma_i(0) =
x_i$ and $\gamma_i(1) = y_{\sigma (i)}$ for some permutation $\sigma
\in\Gamma$. \ere

We are now in a position to prove Theorem \ref{mainstat}. We consider first the case when $\Gamma$ is a transitive subgroup of $\sn$.
Write $\pi : X^n\lrar\gp{n}X$ the projection.  We will proceed in two steps: (i) $\pi_1(\gp{n}X)$ is abelian if $\Gamma$ is transitive and (ii) $\pi_1(\gp{n}X)=H_1(X,\bbz )$.

To prove (i), we will use the trick of Smith which consists in ``lifting and commuting" loops. First of all notice that by Corollary \ref{fixedpoint2}, $\pi$ is surjective
 on fundamental groups. Write
$\eta_{j}, 1\leq j\leq k$ the $j$-th generator of $\pi_1(X)$, and write a
presentation of the product as
$$\pi_1(X)^n = \langle \eta_{11},\ldots, \eta_{1k}, \cdots ,
\eta_{n1},\ldots, \eta_{nk}\ |\ \hbox{relations} \rangle $$ where
$\eta_{ij}$ is the $j$-th generator of the $i$-th copy of $\pi_1(X)$
in $\pi_1(X)^n$.
Let $\alpha$ and $\beta$ be two homotopy classes in
$\pi_1(\gp{n}(X))$. We will show that they commute.
Chooses $\tilde \alpha$ and $\tilde \beta$ in
$\pi_1(X^n) = \pi_1(X)^n$ which project respectively onto $\alpha$ and
$\beta$. Write these two classes additively as
$$\tilde\alpha=\sum_{i,j}{a_{i,j}\eta_{i,j}} \ \ \ ,\ \ \
\tilde\beta=\sum_{i,j}{b_{i,j}\eta_{i,j}}$$
Observe that since $\pi_*(\eta_{i,j}) = \pi_*(g\eta_{i,j})=\pi_*(\eta_{g(i),j})$, $g\in\Gamma$, the
transitivity of the action garantees that
$$\pi_*(\eta_{r, j}) = \pi_*(\eta_{s,j })\ \ \ ,\ \ \ 1\leq
r,s \leq n$$
If we define
$$\tilde\alpha_1 =\sum_{i,j}{a_{i,j}\eta_{1,j}} \ \ \ ,\ \ \
\tilde\beta_2 =\sum_{i,j}{b_{i,j}\eta_{2,j}}$$
then
$$\pi_*(\tilde\alpha_1) = \pi_*(\tilde\alpha) = \alpha
\ \ ,\ \ \pi_*(\tilde\beta_2) = \pi_*(\tilde\beta) = \beta$$
The point of this manipulation is that now $\tilde\alpha_1$ and
$\tilde\beta_2$ commute in $\pi_1(X^n)$ since they are made out of
generators that live in different factors of $\pi_1(X)^n$. Since
these classes commute so do their images $\alpha$ and $\beta$ in
$\pi_1(\gp{n}(X))$. This ``lifting and commuting" argument
proves therefore that $\pi_1(\gp{n}X)$ is abelian.
In fact we have proven a little more

\bco\label{composite} In the case of a transitive action $\Gamma\subset Aut\{1,\ldots, n\}$, the composite
$$\tau : \pi_1(X)\hookrightarrow \pi_1(X)^n\lrar \pi_1(\gp{n}X)=H_1(\gp{n}X)$$
is surjective, where the first map is inclusion as the first factor.
\eco

\begin{proof} This is restating what we previously constructed.
If we refer to the first factor of $\pi_1(X)^n$ by
$\pi_1(X)$ which is generated by the $\eta_{1,j}, 1\leq j\leq k$,
then for $\alpha\in\pi_1(\gp{n}X)$, with preimage $\tilde\alpha\in
\pi_1(X)^n$ and $\tilde\alpha = \sum a_{i,j}\eta_{i,j}$, the class
$\tilde\alpha_1 := \sum a_{i,j}\eta_{1,j}$ is also a preimage.
\end{proof}

Finally we verify that $\pi_1(\gp{n}X) = H_1(\gp{n}X)$ coincides with $H_1(X;\bbz )$ if
$n\geq 2$. The composite $\tau$ in corollary \ref{composite}
factors necessarily through the abelianization
$\pi_*: H_1(X;\bbz)\lrar H_1(\gp{n}X)$ and this map as we indicated is
surjective. It then remains to show it is also injective and this can be
done through a neat little argument.

\ble\label{abelian1} Let $i: X\hookrightarrow\gp{n}X$ be the basepoint inclusion
$x\mapsto [x,*,\ldots, *]$, $n\geq 2$. Then for all $k\geq 0$, $i_* :
H_k(X;\bbz)\lrar H_k(\gp{n}X)$ is a monomorphism. \ele

\begin{proof}
  We show that on cohomology, the induced map $i^*: H^k(\gp{n}X)\lrar H^k(X)$
  is surjective. Let $\alpha\in H^k(X;\bbz )$. Since $X$ is of the homotopy type of a
  CW-complex, the class $\alpha$ can be represented by
  the (homotopy class) of a map $f: X\lrar K(\bbz, k)$; that is $\alpha = f^*(\iota )$ where $\iota\in
  H^k(K(\bbz,k))\cong\bbz$ is a generator. A model for the Eilenberg-MacLane space $K(\bbz, k)$ is $\spy (S^k)$ and
  this is an abelian topological monoid (which we write multiplicatively). We
  can then extend $f$ to a map $\hat f: \gp{n}X\lrar K(\bbz, k)$ sending
  $[x_1,\ldots, x_n]\longmapsto f(x_1)\cdots f(x_n)$. The diagram
$$\xymatrix{
X\ar[r]^f\ar[d]^i&K(\bbz, k)\\
\gp{n}X\ar[ru]_{\hat f}
}
$$ defines a class $\hat\alpha := \hat f^*(\iota)$ such that
$i^*(\hat\alpha ) = \alpha$. Since $i^*$ is surjective, then $i_*$ is
injective.
\end{proof}

The next lemma (and its proof) shows what happens when we go from a single
orbit to multiple ones.

\ble\label{firsthomology} Let $\Gamma\subset\hbox{Aut}\{1,\ldots, n\}$ and
$m$ the total number of orbits under this action. Then
$$H_1(\gp{n}X;\bbz )\cong \left(H_1(X;\bbz)^{\oplus n}\right)^\Gamma
\cong H_1(X;\bbz)^m$$
where the middle term is the submodule of invariants under the action of $\Gamma$ permuting the factors.
\ele

\begin{proof}
We break the action of $\Gamma$ on $\Omega=\{1,\ldots,n\}$ into orbits $\Omega_1,\ldots, \Omega_m$ of respective sizes $k_i\geq 1$. We can change the action of $\Gamma$ up to conjugation so that those
orbits are $\Omega_1=\{1,\ldots, k_1\}$, $\Omega_2=\{k_1+1,\ldots, k_1+k_2\}$,
$\cdots, \Omega_m=\{k_1+\cdots +k_{m-1}+1,\ldots, \sum k_i\}$.  Set $Y_i:=\hbox{Aut}(\Omega_i)
\cong\mathfrak S_{k_i}$. Then the $Y_i$ commute among each other and their product in $\sn$ is a subgroup $Y_1\cdots Y_m$ which is isomorphic to the Young subgroup $\mathfrak S_{k_1}\times\ldots \times\mathfrak S_{k_m}$ (notation as discussed at the start of section \ref{structure}). We can assume from now on and wlog  that
 $\Gamma\subset \mathfrak S_{k_1}\times\ldots \times\mathfrak S_{k_m}$ so that
 $\Gamma$ permutes $\{1,\ldots, k_1\}$, then $\{k_1+1,\ldots, k_1+k_2\}$, etc.
 We decompose $X^{n}= X^{k_1}\times X^{k_2}\times\cdots\times X^{k_m}$
 so that we have the composite
$$X^m\fract{i_m}{\lrar}  X^{n}\fract{\pi}{\lrar} \gp{n}(X)
\fract{q_m}{\lrar} \sp{k_1}(X)\times\cdots\times\sp{k_m}(X)$$
where $\pi$ and $q_m$ are quotient maps and $i_m$ is the map that identifies the $i$-th factor of $X^m$ with the first factor of $X$ in $X^{k_i}$. It is immediate to see that this composite induces an isomorphism on $H_1$ (because each composite $X\lrar X^{k_i}\lrar\sp{k_i}(X)$ does by Lemma \ref{composite}) and so $H_1(\gp{n} ;\bbz)$ contains $H_1(X)^{\oplus m}$ as a subgroup.
In fact we show this is an isomorphism. Indeed we know that
$\pi_*: H_1(X^n)\lrar H_1(\gp{n}(X))$ is surjective by Corollary \ref{fixedpoint2} and by the fact that the fundamental group of $\gp{n}(X)$ is abelian. Given $\alpha\in H_1(\gp{n}(X))$, we choose a lift
$\tilde\alpha = \sum a_{ij}\eta_{ij}\in H_1(X^n)\cong H_1(X)^{\oplus n}$,
where $\eta_{ij}$ is the $j$-th generator in the $i$-th copy of $H_1(X)\hookrightarrow
H_1(X)^{\oplus n}$. Suppose $i$ is in the $t$-th orbit $\Omega_t$, then since the
action of $\Gamma$ on $\Omega_t$ is transitive, we have that
\begin{equation}\label{lift2}
\pi_* (\eta_{ij}) = \pi_*(\eta_{k_1+\cdots +k_{t-1}+1,j})
\end{equation}
This means that a lift of $\alpha$ can be chosen to be in the form
$\sum b_{ij}\eta_{k_1+\cdots + k_i+1,j}$ and hence is in the image of $i_m: H_1(X)^{\oplus m}\hookrightarrow H_1(X^n)$.  The composite
$\pi\circ i_m: H_1(X)^{\oplus m}\lrar H_1(\gp{n}X)$ is therefore surjective and hence $H_1(\gp{n}X)\cong H_1(X)^{\oplus m}$. Remains to see that this is also
isomorphic to $H_1(X^n)^\Gamma\cong \left(H_1(X)^{\oplus n}\right)^\Gamma$ but the
arguments are the same and we leave them to the reader.
\end{proof}

\begin{proof} \textit{of Theorem \ref{mainstat} }:
As in the proof of Lemma \ref{firsthomology} previous, we suppose
$\Gamma$ acts with non-trivial orbits
$\Omega_1,\ldots, \Omega_r$, of lengths $k_1,\ldots, k_r\geq 2$, and the
remaining orbits $n-\sum k_i$ are trivial. We can also assume
that $\Omega_1=\{1,\ldots, k_1\}$, $\Omega_2=\{k_1+1,\ldots, k_1+k_2\}$, etc. The quotient of $X^n$ by $\Gamma$ splits  as $\left(X^{\sum_1^r k_i}/\Gamma\right)\times X^{n-\sum_1^r k_i}$ and we therefore have to determine the fundamental group of the first factor.
By lifting homotopy classes as in Lemma \ref{firsthomology} (see (\ref{lift2})), and
using the ``lifting and commuting" argument explained earlier,
 we show that $\pi_1\left(X^{\sum k_i}/\Gamma\right)$ is abelian.
By Lemma \ref{firsthomology} we can then write
$$\pi_1\left(X^{\sum k_i}/\Gamma\right)\cong
H_1(X^{\sum_1^r k_i}/\Gamma ;\bbz)\cong H_1(X;\bbz)^r$$
and the theorem follows.
\end{proof}

\subsection{Higher Homotopy Groups} We can ask generally how
$\pi_i(\gp{n}X)$ compares to $H_i(X,\bbz )$ for
finite $n, i>1$. We say a simply connected $X$ is $r$-connected if its
homology vanishes up to and including degree $r$.
When $\Gamma = \sn$ and
using simplicial techniques Dold and Puppe show the following \cite{dp}.

\bth\label{stable} Let $X$ be $r$-connected, $r\geq 1$, then for $n>1$
$$\pi_i(\sp{n}X)\cong \tilde H_i(X;\bbz )$$
provided that $0\leq i\leq r+2n-1$.
\end{theorem}

We give a short and novel proof of this theorem based on results in \cite{braids}.

\begin{proof} If $X$ is simply connected, then so are its symmetric products.
It follows that $\pi_1(\sp{n+1}X,\sp{n}X) = 0$. The relative
Hurewicz theorem (see \cite{hatcher}, \S4.2) then says that the first non-zero
$\pi_*(\sp{n+1}X,\sp{n}X)$ is isomorphic to the first non-zero $H_*(\sp{n+1}X,
\sp{n}X)$. Since $H_i(\sp{n+1}X,\sp{n}X)\cong\tilde H_i(\bsp{n+1}X)$ is
zero for $i\leq r+2n$ by Theorem  1.3 of \cite{braids}, it follows that
$\pi_i(\sp{n+1}X,\sp{n}X)$ is zero within the same range. From the long exact
sequence of the pair $(\sp{n}X,\sp{n-1}X)$, we find that
$$\pi_i(\sp{n}X)\fract{\cong}{\lrar} \pi_i(\sp{n+1}X)\ \ \
\hbox{for}\ \ \ i \leq r+2n-1$$
and hence by induction
$\pi_i(\sp{n}X)\cong \pi_i(\spy X)\cong \tilde H_i(X,\bbz)$ for $0\leq i\leq r+2n-1$.
\end{proof}

\bre Theorem \ref{stable} is no longer true if $X$ is not simply connected.
It can be shown for example that for $C$ a closed topological surface of genus $g=2$
there is a Laurent polynomial description\footnote{see
\textit{Remarks on symmetric products of curves}, by the first author,
ArXiv math.AT/0402267.}
$$\pi_2(\sp{2}C) =
\bbz [t_i, t_i^{-1}], \ 1\leq i\leq 4$$
\ere


\section{Bounded Multiplicity Configurations: Constructions and Examples}
\label{configs}

There are naturally defined ${\mathfrak S}_n$-invariant subspaces in
$X^n$ (for $n\geq 2$ a fixed integer) which we list below.  Throughout
$\ast\in X$ denotes the basepoint.
\begin{itemize}
\item The $d$-th fat diagonal
$F_d(X,n)$ is the subspace of $X^n$ of all tuples
$(x_1,\ldots, x_n)$ such that $x_{i_1}=x_{i_2}=\cdots = x_{i_d}$ for
some $d$-set of indices $1\leq i_1 < i_2<\cdots <i_d\leq n$. These are
the ``$d$-equal" subspaces considered for example in \cite{bjorner} in the case
$X=\bbr^2$.
\item
$F^d(X,n) = \{(x_1,\ldots, x_n)\in X^n\ |\ \hbox{no entry $x_i$ can
appear more than $d$-times}\}$. We have
$$F^d(X,n) = X^n - F_{d+1}(X,n)$$
\item The symmetric group $\mathfrak S_n$ acts by permuting coordinates
and preserves these subspaces, thus we can define the unordered analogs of the above
constructions
$$B_d(X,n) = F_d(X,n)/{\mathfrak S}_n\ \ \hbox{and}\ \
B^d(X,n) = F^d(X,n)/{\mathfrak S}_n$$
Clearly $B_{n}(X,n) = X$ while $B_2(X,n)$ is the image of the ``fat
diagonal" in $\sp{n}X$. These are the main objects of study in this paper.
A widespread notation in the literature is to write
$$F(X,n) := F^1(X,n)\ \ \ \ ,\ \ \ B(X,n) = B^1(X,n) $$
for the configuration spaces of ordered (resp. unordered) pairwise distinct points of $X$.
\end{itemize}

Note that $B^d(X,n)$ is not a functorial construction since the constant map
$X\lrar X$, $x\mapsto *$, doesn't induce a self-morphism of $B^d(X,n)$. On the
other hand the fat diagonal $B_d(X,n)$ is functorial and an
invariant of homotopy type.

\bex There is a homeomorphism
$$
B_d(X,n) = X\times\sp{n-d}(X)\ \ \ \ \ \
\hbox{when $d > [n/2]$ and $n\geq 3$}$$
This is true because there can be only one entry with
multiplicity $\left[{n\over 2}\right]+1$ or higher so it can be singled out. \eex

\bex\label{case2} Suppose $[n/3] < d \leq [n/2]$. Then there is a
pushout diagram
\begin{equation}\label{pushout2}
\xymatrix{ (X\times X)\times\sp{n-2d}(X)\ar[r]^{\ \ \ \ \ f}\ar[d]^{\pi\times
    1}&X\times\sp{n-d}(X)\ar[d]^\beta\\ \sp{2}(X)\times\sp{n-2d}(X)
  \ar[r]^{\ \ \ \ \ \alpha}&B_d(X,n) }
\end{equation}
where $f (x,y,z) = (x,y^dz)$, $\beta (a,b) = a^db$ and
$\alpha(xy,z)= (x^dy^dz )$. Here for ease we use the abelian product $x^ry^s$
to denote the point $[x,...x,y...,y]$ ($x$ repeating $r$-times and $y$ repeating $s$ times). The justification behind this pushout is that
when $[n/3] < d \leq [n/2]$, there are at most two points
having multiplicity at least $d$ and when both occur they are
indistinguishable.\eex

The following has been discussed in section \ref{structure}.
We recall that ``Young" refers to a subgroup of $\sn$ conjugate to a standard
embedding of $\mathfrak S_{d_1}\times\ldots \times \mathfrak
S_{d_k}$ in $\sn$, with $\sum d_i= n$. Write $\mathfrak
S_d\subset\sn$ the standard inclusion consisting of the subgroup permuting the
first $d$ letters in a set of $n$ letters, $d\leq n$, leaving the others fixed.  Then

\ble\label{functorial} We have the following identification\
$F_d(X,n) = \overline{X^n}_{(\mathfrak S_d)}
$
\ele

Our next statement is a general observation that fits in naturally in this
section. We shall give a characterization in terms of Young strata $X^n_{({\small
Young})}$ of some
naturally occurring constructions. Below $n\geq 1$ is a fixed integer and all
spaces are path-connected. A space $X$ is called ``locally contractible" if
any $x\in X$ has a neighborhood deformation retract.

\bpr\label{functorialdescription}
 For each locally contractible $X$, suppose we can associate a
$\sn$-invariant subspace
$F(X)\subset X^n$ so that any map $f:X\lrar Y$ induces maps
\begin{equation}\label{functoriality}
\xymatrix{F(X)\ar[d]\ar@{^(->}[r]&X^n\ar[d]^{f^n}\\
F(Y)\ar@{^(->}[r]&Y^n
}
\end{equation}
Then $F(X)$ must be a union of skeleta of the form
$\overline{X^n}_{(H)}$ for some choices of Young-subgroups $H\subset\sn$.
\epr

\begin{proof}
  What the Proposition says for $n=1$ is that $F(X)=X$, and for $n=2$ there are
  two possibilities; either $F(X)=X^2$ or $F(X)$ is the diagonal subspace.
  Note that by choosing $f : \{x\}\lrar X, x\mapsto x$, in diagram
  (\ref{functoriality}), we see that $(x,\ldots, x)\in F(X)$, and $F(X)$
  always contains the diagonal.

We say that a space $X$ is ``$k$-homogeneous" if $\forall (x_1,\ldots, x_k)\in
X^k$, $x_i\neq x_j, i\neq j$, and $\forall (y_1,\ldots, y_k)$, there exists a
\textit{continuous\footnote{Customarily the definition of homogeneity involves
    the existence of homeomorphisms. This is too strong for our purpose.}} map
$h:X\rightarrow X$ such that $h(x_i)=y_i$.  We show that any locally
contractible space is $k$-homogeneous for any $k$.  Pick $(x_1,\ldots, x_k)$
and $(y_1,\ldots, y_k)$ as above.  Remember that a pair of spaces $(X,A)$ has
the ``homotopy extension property" or HEP if for every map $f:X\lrar Y$ and
every homotopy $h: A\times I\lrar Y$ with $h(-,0)=f_{|A}$, there is a homotopy
$H: X\times I\lrar Y$ extending $h$. An NDR pair has the HEP, in particular if
$A$ is a finite set of points in a locally contractible $X$, then $(X,A)$ has
the HEP. Set $A=\{x_1,\ldots, x_k\}\subset X$ and consider the constant map
$f:X\lrar X$ at $x_0\in X$. The restriction of $f$ to $A$ is homotopic to the
map $h_1: A\lrar X$ sending $x_i\mapsto y_i$ (this is because $X$ is
path-connected and any choice of paths from $y_i$ to $x_0$ gives a null
homotopy for $h_1$). the homotopy $h_t: A\lrar X$ from $h_0=f$ to
$h_1$ extends to $H_t : X\lrar X$ and so the
map $H_1 : X\lrar X$ is an extension of $h_1$ with the desired properties.

Going back to the proof, let us pick an element
$(x_1,\ldots ,x_n)\in F(X)$. It must belong to some stratum
$X^n_{(H)}$, $H\cong \mathfrak S_{d_1}\times\ldots \times \mathfrak S_{d_k}$, since
the Young strata form a stratification. The
transitivity property established above shows that there are selfmaps
of $X$ that send $(x_1,\ldots, x_n)$ to any tuple $(y_1,\ldots, y_n)$
provided that $y_{i_1}=\ldots =y_{i_d}$ if $x_{i_1}=\cdots =x_{i_d}$.
This immediately implies that not only the whole stratum is in $F(X)$, but
also its closure $\overline{X^n}_{(H)}$. From there we deduce that
$F(X)$ must be a union of strata $\overline{X^n}_{(H)}$ as asserted.
\end{proof}


\subsection{Cell and Manifold Structures}\label{cellandmanifold}

When $X$ is a simplicial complex (i.e. a triangulated space), then $\gp{n}X$
has the structure of a CW complex.  This cell structure for the case of
symmetric products is discussed in \cite{dold,liao} for example.  These
descriptions lead to a more precise statement.

\bpr\label{celldecomp} For a simplicial complex $X$, $\gp{n}X$ is a CW
complex such that the skeleta of its associated orbit
stratification are subcomplexes.  \epr

The point is the existence of a $\sn$-invariant (thus $\Gamma$-invariant)
simplicial structure on $X^n$, induced from one on $X$, such that the fat
diagonal is a CW-subcomplex. Because this simplicial structure is invariant,
it induces a \textit{cellular} (not necessarily simplicial) structure on the
quotient $\gp{n}X$ and such that the quotient map $\pi :
X^n\rightarrow\gp{n}X$ is cellular.

We recall a space $X$ is $r$-connected, $r\geq 0$, if the homotopy groups
$\pi_i(X)$ vanish for $i\leq r$.

\bco\label{connectivity} If $X$ is an $r$-connected simplicial
complex, then $B_d(X,n)$ is also $r$-connected. \eco

\begin{proof} The proof of (\cite{ks2}, corollary 3.5) applies verbatim to this
  situation. It uses the CW decomposition on $B_d(X,n)$ as a subcomplex of
  $\sp{n}X$ (Proposition \ref{celldecomp}).
\end{proof}

Even if $X$ is a manifold, the spaces $B_d(X,n)$ are rarely manifolds
themselves. If $\dim X=1$, then two cases occur: $X=I$ is an open interval
or $X=S^1$. When $X$ is homeomorphic to an open interval
$\stackrel{\circ}{I}=]0,1[$, then $B_d(\stackrel{\circ}{I},n)$ is the obvious subspace of the
$n$-simplex
$$\sp{n}(\stackrel{\circ}{I}) = \{(t_1,\ldots, t_n), 0<t_1\leq\cdots\leq t_n <1\}
=\Delta_n$$
consisting of part of the boundary where $t_{i_1}=\cdots =t_{i_d}$
for some subset $\{i_1,\ldots, i_d\}$.
When $X=S^1$, then by a refined version of the result of Morton (Remark \ref{morton}),
$B_d(S^1,n)$ is a bundle over $S^1$ with fiber some faces of $\Delta_{n-1}$. For certain values of $d$ this is a manifold. For example
$B_2(S^1,n)\subset\sp{n}(S^1)$ is an $n-2$-sphere bundle over $S^1$
and this is a manifold.

When $\dim X\geq 2$, the following theorem gives complete conditions on when $B_d(X,n)$
are manifolds.

\bth\label{manifold}
Let $X$ be a  manifold of dimension $m\geq 2$ and suppose $n\geq 3$ , $d\geq 2$.
Then $B_d(X,n)$ is a
 manifold again if and only if one of the following conditions is
satisfied:\\
(i) $d=n$ or $d=n-1$, \\
(ii) $d> [n/2]$ and $X$ is a topological surface.
\end{theorem}

\begin{proof} The ``if'' part is clear since (i) $B_n(X,n)=X,
B_{n-1}(X,n)\cong X\times X$, and (ii) $B_d(X,n)\cong X\times\sp{n-d}(X)$ when $d> [n/2]$,
and the symmetric products of topological surfaces are manifolds.
We need to see that aside from these special cases, $B_d(X,n)$ cannot be a
manifold. Let's assume at the start that $m\geq 2$ and
that $d\geq 3$. Here we can set $n\geq d+2$ (the cases $n=d$ and
$n=d+1$ having been settled).
Pick a point $\zeta = [x,\cdots ,x, y,y, z_1,\ldots, z_{n-d-2}]\in B_d(X,n)$, $x$ repeating
$d$-times. We assume that $x\in V$, $y\in U$ and
$z_i\in U_i$ all lie in neighborhoods $V,U$ and $U_i$ that are pairwise disjoint and
homeomorphic to $\bbr^m$, $m=\dim X$.
Then an open neighborhood of
$\zeta\in B_d(X,n)$ is
$V\times\sp{2}(U)\times U_1\times\cdots \times U_{n-d-2}$. This
can be given a Euclidean structure if and only if $\sp{2}(U)\cong
\sp{2}(\bbr^m)$ can be given a Euclidean structure. Unless $m=1$ or $m=2$, it
is not possible to find a neighborhood of the origin $[0,0]\in\sp{2}(\bbr^m)$
that is homeomorphic to $\bbr^{2m}$.
The argument here was already given at the start of \S\ref{diagonalstratum} and is based
on showing that $\sp{2}(\bbr^m)\cong \bbr^m\times c\bbr P^{m-1}$.
To see this one decomposes $W=\bbr^m\oplus\bbr^m$ as $D\oplus V$ where
$D$ is the diagonal subspace fixed by the transposition $A(x,y)=(y,x)$ acting
on $W$, and where $V=\{(x,-x), x\in\bbr^m\}$ is the orthogonal subspace. Now $V$ is a copy of $\bbr^m$ on which $A$ acts via $x\mapsto -x$. The quotient
of $D\oplus V$ by this action $\bbz_2$ (generated by $A$)
is clearly $\bbr^m\times c\bbr P^{m-1}$.
Since no neighborhood of the cone point $[0,0]$ can
be homeomorphic to $\bbr^{2m}$ if $m\geq 3$,
the only possibility is therefore $m=2$ (or $m=1$ but this is not relevant as $m\geq 2$)
and $X$ is a topological surface.

We are then left to discuss the case $d=2, m\geq 2$.
Assume without loss of generality that $n=2d=4$.
Let $x, y$ be two points of $X$ lying in disjoint Euclidean neighborhoods $U$ and $V$
respectively. A neighborhood of $x^2y^2\in B_2(X,4)$ can be seen to be the pushout
\begin{equation}\label{voisinage}
U\times\sp{2}(V)\cup_{U\times V}\sp{2}(U)\times V
\end{equation}
where $U\times V$ includes in $B_2(X,4)$ via the map $(x,y)\mapsto
x^2y^2$. This neighborhood has an embedded copy of $U\times V\cong
\bbr^m\times\bbr^m$. This embedded copy is of codimension $m\geq 2$.
Therefore if the neighborhood in (\ref{voisinage}) were homeomorphic to Euclidean space,
removing $U\times V$ from this neighborhood
wouldn't disconnect it. But this isn't the case since the complement of
$U\times V$ in this pushout is disconnected.
\end{proof}


\subsection{Homology}
We make some useful remarks about the homological structure of $B_d(X,n)$
relevant to \S\ref{fundgroupbd}.
First of all and using the known fact that $H_*(\sp{n}X;\bbz)$ embeds in $H_*(\sp{n+d}X;\bbz )$
for $d\geq 1$ via the ``basepoint embedding"
$\sp{n}(X)\hookrightarrow\sp{n+d}(X)$ sending
$[x_1,\ldots, x_n]\longmapsto [x_1,\ldots, x_n,\ast ,\ldots, \ast ]$,
where $*$ is a basepoint, we immediately deduce the following

\ble\label{embeds} The inclusion $\sp{n}X\hookrightarrow B_d(X,n+d)$,
$[x_1,\ldots, x_n]\mapsto [x_1,\ldots, x_n,\ast ,\ldots, \ast]$ induces a
monomorphism in homology. \ele

Generally $H_*(B_d(X,n))$ doesn't embed in $H_*(B_d(X,n+1))$. This
is illustrated in the case $d=2$ by the following calculation for example.
Here we recall that $B_2(S^n,3)\cong S^n\times S^n$, while

\bpr\label{walidt}\cite{taamallah} $B_2(S^{n},4)$ has the integral homology of
$$S^n\vee
\Sigma^{n+1}\bbr P^{n-1}\vee
\Sigma^{n+1}\bbr P^{n-1}
\vee \Sigma^{2n+1}\bbr P^{n-1}\vee \Sigma^{2n-1}
\bbr P^2$$
\epr

\begin{proof} There is an
identification due to Steenrod of the mapping cone
of the basepoint inclusion $S^n\hookrightarrow\sp{2}(S^n)$
with the suspension $\Sigma^{n+1}\bbr P^{n-1}$ (see \cite{hatcher}, example 4K.5).
This gives that $H_*(\sp{2}(S^n))$ has a torsion free class in degree $n$
and another in degree $2n$ only if $n$ is even,
and a copy of $\bbz_2$ in degrees $r=n+2k$ for $1\leq k\leq [(n-1)/2]$. The rest of
the argument follows from a careful analysis of the Mayer-Vietoris exact sequence
associated to the pushout diagram in Example \ref{case2}. Details in \cite{taamallah}.
\end{proof}

Note that $H_*(B_2(S^n,4))$ doesn't satisfy Poincar\'e duality consistently
with Theorem \ref{manifold}. The following more general calculation for fat diagonals
of spheres is valid with rational coefficients.

\bpr\label{ft}\cite{ft} The rational homotopy type of the fat diagonal of an
odd sphere is
$$B_2(S^{2k+1},n) \simeq_\bbq\begin{cases}S^{2k+1},& n\ \hbox{even}\\
S^{2k+1}\times S^{m(2k+1)-1},& n=2m+1\ \hbox{odd}
\end{cases}
$$
\epr

Note how this is naturally consistent with Proposition \ref{walidt} when $n=4$ since
$\bbr P^{2k}$ has the rational homology of a point.

\section{Fat diagonals and their fundamental groups}\label{fundgroupbd}

In this last section we prove Theorem \ref{main1}.  We will be using the
Van-Kampen Theorem in the following form: let $G =
\langle g_1,\ldots, g_n\ |\ r_G\rangle $ and $H=\langle h_1,\ldots, h_m\ |\
r_H\rangle $ be finitely presented groups, where $r_G$ and $r_H$ stand for the
corresponding ideals of relations. Let $K$ have presentation $K = \langle
k_1,\ldots, k_t\ |\ r_K\rangle $, mapping to $H$ and $G$ via $\alpha$ and $\beta$
respectively.  The amalgamated product $H*_KG$ has presentation
$$\langle g_1,\ldots, g_n, h_1,\ldots ,h_m\ |\ r_G, r_H,
\alpha (k_i)\beta (k_i)^{-1}\rangle $$
Note that the relations $r_K$ do not matter. Note also that if $\alpha$ is
surjective, then the ``universal" map $G\lrar H*_KG$ is surjective as well.

Here is now the main assertion of this section

\ble\label{bdsecond} For $n\geq 2$ and $1\leq d\leq {n\over 2}$,
$\pi_1(B_d(X,n))$ is abelian.  \ele

\begin{proof} We can assume $d\geq 2$ since when $d=1$, $B_1(X,n)=\sp{n}X$.
  As in the proof of the commutativity of $\pi_1(\gp{n}X)$ in
  \S\ref{fundgroup}, the idea here is to start with two elements in
  $\pi_1(B_d(X,n))$, lift them to $\pi_1(F_d(X,n))$ via the surjection
  $\pi_*:\pi_1(F_d(X,n))\twoheadrightarrow \pi_1(B_d(X,n))$ (see Corollary \ref{fixedpoint2}), then using the
  action of $\sn$ separate them into elements that commute.  Write
$$\pi_1(X)=G = \langle g_i, i\in J\ |\ \hbox{relations}\rangle \ \ ,
\ \hbox{with $J$ some finite
indexing set}$$ For
$I = \{i_1,\ldots, i_d\}$, define
$$\Delta_I(X,n)=\Delta_{i_1, \cdots, i_d}(X,n) := \{(x_1,\ldots, x_n)\ |\
x_{i_1}=\cdots = x_{i_d}\}$$
Then
$\pi_1(\Delta_{i_1, \cdots, i_d}(X,n))\subset G^n$ has generators
$$g^I_i\ ,\ g^I_{a,j}\ \ \ ,\ \ i,j\in J, a\not\in I$$ where $g^I_i =
(e,\cdots ,g_i,\cdots , g_i,\cdots , e)$ is the element of $G^n$ with $g_i$ in
the $i_1,i_2,\cdots ,i_d$ entries and the identity $e$ in the remaining
entries, while $g^I_{a,j}=(e,\cdots ,g_j,e ,\cdots , e)$ is the element of
$G^n$ with $g_j$ in the $a^{th}$ entry and the identity $e$ in the remaining
entries.  Under the isomorphism $\pi_1(\Delta_{i_1, \cdots, i_d}(X,n))\cong
G\times G^{n-d}$, the $g^I_i$ generate the first factor $G$ and the
$g^I_{a,j}$ generate the second factor $G^{n-d}$. For example
$\pi_1(\Delta_{1,2}(X,3))$ has generators
$$g^{12}_i\ ,\ g^{12}_{3,j}\ \ ,\ i,j\in J$$

Since $F_d(X,n)=\bigcup\Delta_I(X,n)$, with $I$ a sequence of the form $1\leq
i_1<\cdots <i_d\leq n$, an iterated use of the Van-Kampen theorem shows that
$\pi_1(F_d(X,n))$ is generated by the $g^I_i$ and $g^{I}_{a,j}$ with $I$
varying over all such sequences, $a\not\in I$ and $i,j\in J$, subject to
various relations. Among these relations are those coming from
$\pi_1(\Delta_I(X,n))$ which are:
\begin{itemize}
\item $g^I_i$ and $g^I_{a,j}$, $a\not\in I$,
commute since they come from different factors of $G^n$.
\item Similarly $g^I_{a,r}$ and $g^I_{b,s}$ commute if $a\neq b$.
\end{itemize}

Write $I_d = \{1,2,\ldots, d\}$ and $I^d = \{n-d+1, n-d+2,\ldots, n\}$.  These
are disjoint subsets since $n\geq 2d$. Consider the pushout diagram
$$\xymatrix{
  X\times X\times X^{n-2d}\ar[r]^f\ar[d]^g&\Delta_{I_d}(X,n)\ar[d]\\
  \Delta_{I^d}(X,n)\ar[r]&\Delta_{I^d}(X,n)\cup \Delta_{I_d}(X,n) }$$ where
the image of $(x, y, (z_1,\ldots, z_{n-2d}))$ under both inclusions $f$ and
$g$ is
$$(x,x,\ldots, x,z_1,\ldots, z_{n-2d},y,y,\ldots, y)$$
If we write $\pi_1(X\times X\times X^{n-2d})$ as
$G\times G\times G^{n-2d}$, then the generators $g_{1,i}$ coming from the
first copy of $G$ map to
$$f_*(g_{1,i}) = g^{I_d}_i\ \ \ \ ,\ \ \ \ g_*(g_{1,i}) = g^{I^d}_{1,i}+
g^{I^d}_{2,i}+\cdots + g^{I^d}_{d,i}$$
By the Van-Kampen theorem, this implies that in $\pi_1(F_d(X,n))$,
the generators of the same name $g^{I_d}_i$ and $g^{I^d}_{a,i}$
satisfy the relation for given $i\in J$
\begin{equation}\label{relation}
g^{I_d}_i = g^{I^d}_{1,i}+g^{I^d}_{2,i}+\cdots + g^{I^d}_{d,i}
\end{equation}

Let now $\alpha, \beta$ be two elements in $\pi_1(B_d(X,n))$ and write
$\tilde\alpha, \tilde\beta$ in $\pi_1(F_d(X,n))$ such that
$\pi_*(\tilde\alpha ) = \alpha$ and $\pi_*(\tilde\beta )=\beta$.  In
terms of generators we can express $\tilde\alpha$ as a finite sum
$$\tilde\alpha = \sum n^I_ig^I_i + \sum m^{I'}_jg^{I'}_{a,j}\ \ \ ,\ \ \
n^I,m^{I'}\in\bbz$$
for some choice of $I,I',a,i,j$. We know that under the $\sn$ quotient,
\begin{equation}\label{sigmaaction}
\pi_*(g^I_i) = \pi_*(g^{\sigma (I)}_i)\ \ \ ,\ \ \
\pi_*(g^I_{a,j}) = \pi_*(g^{\sigma (I)}_{\sigma (a),j})\ ,\ a\not\in I
\end{equation}
The first equality follows from the fact that $\sn$ acts $d$-transitively
on $\{1,\ldots, n\}$ meaning that any subset of cardinality $d$ is mapped
to any other by the action. From the relations (\ref{sigmaaction}), another
choice of a lift for $\alpha$ would then be
$$
\tilde\alpha_1 = \sum n_ig^{I^d}_i + \sum m_jg^{I^d}_{1,j}
$$ and this lives in $\pi_1(\Delta_{I^d}(X,n))$.  We replace
$g^{I^d}_i$ by $g^{I_d}_i$ and this in turn can be replaced by
$g^{I^d}_{1,i}+\cdots + g^{I^d}_{d,i}$ according to (\ref{relation}).
Similarly we can replace each $g^{I^d}_{a,i}$ by $g^{I^d}_{1,i}$ since
one is mapped into the other by the transposition $(1,a)$. At the end
and for $n\geq 2d$ there is a choice of lift
\begin{equation}\label{tildealpha}
\tilde\alpha_2 =  \sum t_kg^{I^d}_{1,k}\ \ \ t_k\in\bbz
\end{equation}
such that $\pi_*(\tilde\alpha_2 ) = \alpha$.
Similarly there is a choice of a lift for $\beta$ of the form
\begin{equation}\label{tildebeta}
\tilde\beta_2 =  \sum \ell_kg^{I^d}_{2,k}\ \ \ \ell_k\in\bbz
\end{equation}
This choice is possible since $d\geq 2$.
The expressions in
(\ref{tildealpha}) and (\ref{tildebeta}) commute in $\pi_1(F_d(X,n))$
and hence $\alpha$ commutes with $\beta$ in $\pi_1(B_d(X,n))$. This shows
that $\pi_1(B_d(X,n))$ is abelian.
\end{proof}

The following corollary implies Theorem \ref{main1}.

\bco Assume $2\leq d\leq n/2$.
Then the map $X\hookrightarrow B_d(X,n)$, $x\mapsto [x,\ast,\ast,\cdots, \ast]$
induces an isomorphism on $H_1$.
\eco

\begin{proof} Equation (\ref{tildealpha}) above shows that if
  $X\hookrightarrow\Delta_{I^d}(X,n)$ is the inclusion of the first factor (basepoint everywhere else),
  then the composite $X\lrar\Delta_{I^d}(X,n)\lrar B_d(X,n)$ is surjective on
  $\pi_1$ and thus on
  $H_1$. It remains to see that $H_1(X) \lrar H_1(B_d(X,n))$ is injective but
  this is the content of Lemma \ref{embeds}.
\end{proof}


A subgroup of $\sn$ is ``$d$-transitive" if it takes any \textit{ordered} $d$
points in $\{1,\ldots n\}$ to any other ordered $d$-points. The following
more general result follows from an adaptation of the proof of Theorem \ref{bdsecond}.

\bco\label{dtransitive} If $\Gamma$ is $d$-transitive and $1\leq d\leq n/2$, then
$\pi_1(F_d(X,n)/\Gamma)$ is abelian.
\eco

\vskip 20pt

\bibliographystyle{plain[8pt]}

\end{document}